\documentclass[12pt]{article}
\usepackage{amsfonts}

\usepackage{amsmath,amssymb,amsthm}
\usepackage{multicol}
\usepackage{color}
\usepackage{hyperref}
\usepackage{graphicx}
\usepackage[english]{babel}

\newtheorem{Def}{Definition}[section]
\newtheorem{Cnj}[Def]{Conjecture}
\newtheorem{Prop}[Def]{Property}
\newtheorem{example}{Example}[section]

\begin{document}
\title{Fan, splint and branching rules.}

\author{V.D.~Lyakhovsky$^1$, A.A.~Nazarov$^{1,2}$ \\
  {\small $^1$ Department of High-Energy and Elementary Particle Physics, SPb State University}\\
  {\small 198904, Saint-Petersburg, Russia}\\
  {\small e-mail: lyakh1507@nm.ru}\\
  {\small$^{2}$ Chebyshev Laboratory,}\\
  {\small Department of Mathematics and Mechanics, SPb State University}\\
  {\small 199178, Saint-Petersburg, Russia}\\
  {\small email: antonnaz@gmail.com}}
\maketitle

\begin{abstract}
Splint of root system for simple Lie algebra appears naturally in
studies of (regular) embeddings of reductive subalgebras. Splint can
be used to construct branching rules. We demonstrate that
splint properties implementation drastically simplify calculations of
branching coefficients.
\end{abstract}

\section{Introduction}
\label{sec:Introduction}

Embedding $\phi$ of a root system $\Delta_1$ into a root system
$\Delta$ is a bijective map of roots of $\Delta_{1}$ to a (proper)
subset of $\Delta$ that commutes with vector composition law in
$\Delta_{1}$ and $\Delta$.
\begin{equation*}
\phi:\Delta_1 \longrightarrow \Delta
\end{equation*}
\begin{equation*}
\phi \circ (\alpha + \beta) =\phi \circ \alpha + \phi \circ \beta,
\,\,\, \alpha,\beta \in \Delta_1
\end{equation*}

Note that the image $Im(\phi)$ must not inherit the root system
properties except the addition rules equivalent to the addition
rules in $\Delta_{1}$ (for pre-images). Two embeddings $\phi_1$ and $\phi_2$
can splinter $\Delta$  when the latter can be presented
as a disjoint union of images $Im(\phi_1)$ and $Im(\phi_2)$.
The term {\it splint} was
introduced by D. Richter in  \cite{richter2008splints} where the
classification of splints for simple Lie algebras was obtained.
There was also mentioned that splint must have tight
connections with the injection fan construction. The fan $\Gamma
\subset \Delta$ was
introduced in \cite{lyakhovsky1996rra} as a subset of root system describing recurrent
properties of branching coefficients for maximal embeddings. Injection fan is an
efficient tool to study branching rules. Later this construction
was generalized to non-maximal embeddings and infinite-dimensional
Lie algebras in \cite{2010arXiv1007.0318L, ilyin812pbc}.

In the present paper we study connections between splint and
injection fan for regular embedding of reductive subalgebras
${\mathfrak a}$ in simple Lie algebra $\mathfrak{g}$. We show that (under
certain conditions described in section \ref{sec:stems and
multiplicity functions}) splint is a natural tool to study
reduction properties of ${\mathfrak g}$-modules
with respect to a subalgebra ${\mathfrak a}\longrightarrow
{\mathfrak g}$. Using this tool we obtain the main result -- the one-to-one correspondence between weight
multiplicities in irreducible modules of splint and branching
coefficients for a reduced module $L^{\mu}_{{\mathfrak
g}\downarrow {\mathfrak a}}$.

\section{Injections and splints}

\label{sec:Injections and splints}

Consider a simple Lie algebra $\mathfrak{g}$ and its regular subalgebra $%
\mathfrak{a}\hookrightarrow \mathfrak{g}$ such that $\mathfrak{a}$
is a
reductive subalgebra $\mathfrak{a \subset g}$ with correlated root spaces: $%
\mathfrak{h}_{\mathfrak{a}}^{\ast }\subset \mathfrak{h}_{\mathfrak{g }%
}^{\ast }$. Let $\mathfrak{a}^{\mathfrak{s}}$ be a semisimple summand of
$\mathfrak{a}$,
this means that $\mathfrak{a}=\mathfrak{a}^{\mathfrak{s}} \oplus \mathfrak{u}(1)\oplus %
\mathfrak{u}(1)\oplus \dots$. We shall consider $\mathfrak{a}^{\mathfrak{s}}$
to be a proper regular subalgebra and $\mathfrak{a}$ to be the
maximal subalgebra with $\mathfrak{a}^{\mathfrak{s}}$ fixed that is the rank
$r$ of $\frak{a}$ is equal to that of $\mathfrak{g}$.

The following notations are used:

$r$ , $\left( r_{\mathfrak{a}^{\mathfrak{s}}}\right) $ --- the rank of
$\frak{g}$ $\left( \mathrm{{resp. \,}\mathfrak{a}^{\mathfrak{s}}}\right) $ ;

$\Delta $ $\left( \Delta _{\frak{a}}\right) $--- the root system;
$\Delta ^{+} $ $\left( \mathrm{{resp. \,}\Delta
_{\frak{a}}^{+}}\right) $--- the positive root system (of
$\frak{g}$ and $\frak{a}$ respectively);

$S \, ,\quad \left( S_{\frak{a}}\right) $ --- the system of simple roots (of $%
\frak{g}$ and $\frak{a}$ respectively);

$\alpha _{i}$ , $\left( \alpha _{\left( \frak{a}\right) j}\right) $ --- the $%
i$-th (resp. $j$-th) simple root for $\frak{g}$  $\left( \mathrm{{resp. \,} \frak{%
a}}\right) $; $i=0,\ldots ,r$,\ \ $\left( j=0,\ldots ,r_{\mathfrak{a}%
^S}\right) $;

$\omega _{i}$ , $\left( \omega _{\left( \frak{a}\right) j}\right) $ --- the $%
i$-th (resp. $j$-th) fundamental weight for $\frak{g}$ $\left( \mathrm{{resp. \,}\frak{%
a}}\right) $; $i=0,\ldots ,r$,\ \ $\left( j=0,\ldots ,r_{\mathfrak{a}%
^S}\right) $;

$W$ , $\left( W_{\frak{a}}\right) $--- the corresponding Weyl group;

$C$ , $\left( C_{\frak{a}}\right) $--- the fundamental Weyl chamber;

$\bar{C}, \left(\bar{C_{\frak{a}}}\right)$ --- the closure of the
fundamental Weyl chamber;

$\epsilon \left( w\right) :=\left( -1\right) ^{\mathrm{length}(w)}$;

$\rho $\ , $\left( \rho _{\frak{a}}\right) $\ --- the Weyl vector;

$L^{\mu }$\ $\left( L_{\frak{a}}^{\nu }\right) $\ --- the integrable module
of $\frak{g}$ with the highest weight $\mu $\ ; (resp. integrable $\frak{a}$
-module with the highest weight $\nu $ );

$\mathcal{N}^{\mu }$ , $\left( \mathcal{N}_{\frak{a}}^{\nu }\right) $ ---
the weight diagram of $L^{\mu }$ (resp. ${}L_{\frak{a}}^{\nu }$ );

$P$ (resp. $P_{\frak{a}} $) \ --- the weight lattice;

$P^{+}$ (resp. $P_{\frak{a}}^{+} $) \ --- the dominant weight lattice;

$\mathcal{E}$ (resp. $\mathcal{E}_{\frak{a}} $) \ --- the formal algebra;

$m_{\xi }^{\left( \mu \right) }$ , $\left( m_{\xi }^{\left( \nu \right)
}\right) $ --- the multiplicity of the weight $\xi \in P$ \ $\left( \mathrm{{%
resp. }\in P_{\frak{a}}}\right) $ in the module $L^{\mu }$ , (resp. $\xi \in
L_{\frak{a}}^{\nu } $);

$ch\left( L^{\mu }\right) $ (resp. $\mathrm{ch}\left( L_{\frak{a}}^{\nu
}\right) $)--- the formal character of $L^{\mu }$ (resp. $L_{\frak{a}}^{\nu
} $);

$ch\left( L^{\mu }\right)  =\frac{\sum_{w\in W}\epsilon (w)e^{w\circ (\mu
+\rho )-\rho }} {\prod_{\alpha \in \Delta ^{+}} \left( 1-e^{-\alpha }\right)
} $ --- the Weyl formula;

$R:=\prod_{\alpha \in \Delta ^{+}}\left( 1-e^{-\alpha }\right) \quad $
(resp. $R_{\frak{a}}: =\prod_{\alpha \in \Delta_{\frak{a}}^{+}} \left(
1-e^{-\alpha }\right) $ ) --- the Weyl denominator.

Let $L^{\mu }$ be completely reducible with respect to $\frak{a}$,
\[
L_{\frak{g}\downarrow \frak{a}}^{\mu }=\bigoplus\limits_{\nu \in P_{\frak{a}%
}^{+}}b_{\nu }^{\left( \mu \right) }L_{\frak{a}}^{\nu }.
\]
\begin{equation}
\pi _{\frak{a}}ch\left( L^{\mu }\right) =\sum_{\nu \in P_{\frak{a}%
}^{+}}b_{\nu }^{(\mu )}ch\left( L_{\frak{a}}^{\nu }\right) .
\label{branching1}
\end{equation}
For the modules we are interested in the Weyl formula for $\mathrm{ch}\left(
L^{\mu }\right) $ can be written in terms of singular elements \cite
{humphreys1997introduction}
\[
\Psi ^{\left( \mu \right) }:=\sum\limits_{w\in W}\epsilon (w)e^{w(\mu +\rho
)-\rho },
\]
namely,
\begin{equation}
\mathrm{ch}\left( L^{\mu }\right) =\frac{\Psi ^{\left( \mu \right) }}{\Psi
^{\left( 0\right) }}=\frac{\Psi ^{\left( \mu \right) }}{R}.
\label{Weyl-Kac2}
\end{equation}
The same is true for submodules $\mathrm{ch}\left( L_{\frak{a}}^{\nu
}\right) $ in (\ref{branching1})
\[
\mathrm{ch}\left( L_{\frak{a}}^{\nu }\right) =\frac{\Psi _{\frak{a}}^{\left(
\nu \right) }}{\Psi _{\frak{a}}^{\left( 0\right) }}=\frac{\Psi _{\frak{a}%
}^{\left( \nu \right) }}{R_{\frak{a}}},
\]
with
\[
\Psi _{\frak{a}}^{\left( \nu \right) }:=\sum\limits_{w\in W_{\frak{a}%
}}\epsilon (w)e^{w(\nu +\rho _{_{\frak{a}}})-\rho _{_{\frak{a}}}}.
\]

Applying formula (\ref{Weyl-Kac2}) to the branching rule (\ref{branching1})
we get a relation connecting the singular elements $\Psi ^{\left( \mu
\right) }$ and $\Psi _{\frak{a}}^{\left( \nu \right) }$ :
\begin{eqnarray}
\frac{\sum_{w \in W}\epsilon (w )e^{w (\mu +\rho )-\rho }}{\prod_{\alpha \in
\Delta ^{+}}(1-e^{-\alpha })} &=&\sum_{\nu \in P_{\frak{a}}^{+}}b_{\nu
}^{(\mu )}\frac{\sum_{w \in W_{\frak{a}}}\epsilon (w )e^{w (\nu +\rho _{%
\frak{a}})-\rho _{\frak{a}}}}{\prod_{\beta \in \Delta _{\frak{a}%
}^{+}}(1-e^{-\beta })},  \nonumber  \label{eq:4} \\
\frac{\Psi ^{\left( \mu \right) }}{R} &=&\sum_{\nu \in P_{\frak{a}%
}^{+}}b_{\nu }^{(\mu )}\frac{\Psi _{\frak{a}}^{\left( \nu \right) }}{R_{%
\frak{a}}}.  \label{singular main}
\end{eqnarray}

In \cite{2010arXiv1007.0318L} it was proven that singular branching coefficients
$k_{\xi }^{\left( \mu \right) }$ corresponding to the injection $\frak{a}%
\hookrightarrow \frak{g}$ are subject to the set of recurrent relations:
\begin{equation}
\begin{array}{c}
k_{\xi }^{\left( \mu \right) }=-\frac{1}{s\left( \gamma _{0}\right) }\left(
\sum_{u\in W/W_{\perp}}\epsilon (u)\;\dim \left( L_{\frak{a}_{\perp }}^{\mu _{\frak{a}%
_{\perp }}\left( u\right) }\right) \delta _{\xi -\gamma _{0},\pi _{%
\widetilde{\frak{a}}}(u(\mu +\rho )-\rho )}+\right. \\
\left. +\sum_{\gamma \in \Gamma _{\widetilde{\frak{a}}\rightarrow \frak{g}%
}}s\left( \gamma +\gamma _{0}\right) k_{\xi +\gamma }^{\left( \mu \right)
}\right) .
\end{array}
\label{recurrent rel}
\end{equation}
where $\frak{a}_{\perp }$ is the subalgebra determined by the roots of $%
\frak{g}$ orthogonal to roots of $\frak{a}$ and $W_{\perp}$ is a Weyl group of $\mathfrak{a}_{\perp}$
\begin{eqnarray}
\Delta _{\frak{a}_{\perp }} &:&=\left\{ \beta \in \Delta _{\frak{g}}|\forall
h\in \frak{h}_{\frak{a}};\beta \left( h\right) =0\right\} ,
\label{delta a ort}
\end{eqnarray}
\begin{eqnarray}
\widetilde{\frak{a}_{\perp }} :=\frak{a}_{\perp }\oplus \frak{h}_{\perp }
\qquad \widetilde{\frak{a}} :=\frak{a}\oplus \frak{h}_{\perp }
\end{eqnarray}
and $\pi$ is the projection operator. Inside the main Weyl chamber $C_{\mathfrak{a}}$
singular branching coefficients coincide with branching coefficients:
$ b_{\xi }^{\left( \mu \right) }=k_{\xi }^{\left( \mu \right) }\; \forall \xi\in C_{\mathfrak{a}}$.
When an injection is maximal the
projection becomes trivial and the relation (\ref{recurrent rel}) is
simplified:
\begin{equation}
\begin{array}{c}
k_{\xi }^{\left( \mu \right) }=-\frac{1}{s\left( \gamma _{0}\right) }\left(
\sum_{u\in W}\epsilon (u) \delta _{\xi -\gamma _{0}, u(\mu +\rho )-\rho
}+\right. \\
\left. +\sum_{\gamma \in \Gamma _{\frak{a}\rightarrow \frak{g}}}s\left(
\gamma +\gamma _{0}\right) k_{\xi +\gamma }^{\left( \mu \right) }\right) .
\end{array}
\label{recurrent relation max}
\end{equation}
The recursion is goverened by the set $\Gamma _{\frak{a}\rightarrow \frak{g}}
$ called the injection fan. The latter is defined by the carrier set $%
\left\{ \xi \right\} _{\frak{a}\rightarrow \frak{g}}$ for the coefficient
function $s(\xi )$
\[
\left\{ \xi \right\} _{\frak{a}\rightarrow \frak{g}}:=\left\{ \xi \in P_{%
\frak{a}}|s(\xi )\neq 0\right\}
\]
appearing in the expansion
\begin{equation}
\prod_{\alpha \in \Delta ^{+}\setminus \Delta _{\frak{a}}^{+}}\left( 1-e^{
-\alpha }\right) =-\sum_{\gamma \in P_{\frak{a}}}s(\gamma )e^{-\gamma };\quad
\label{product}
\end{equation}

Now we remind two definitions introduced in \cite{richter2008splints}

\begin{Def}
Suppose $\Delta _{0}$ and $\Delta $ are root systems with corresponding
weight lattices $P_{0}$ and $P$. Then $\phi $ is an ``embedding'',
\begin{equation}
\phi :\left\{
\begin{array}{l}
\Delta _{0}\hookrightarrow \Delta , \\
P_{0}\hookrightarrow P,
\end{array}
\right.
\end{equation}
if \newline
\noindent (a) it injects $\Delta _{0}$ in $\Delta $, and \newline
\noindent (b) acts homomorphically with respect to the vector groups in $%
P_{0}$ and $P$:
\[
\phi (\gamma )=\phi (\alpha )+\phi (\beta )
\]
for any triple $\alpha ,\beta ,\gamma \in P_{0}$ such that $\gamma =\alpha
+\beta $.
\end{Def}

$\phi$ induces an injection of formal algebras $:{\mathcal{E}}_0
\hookrightarrow \mathcal{E}$ and for the image ${\mathcal{E}}%
_i=Im_{\phi}\left( {\mathcal{E}}_0\right)$ one can consider its inverse $%
\phi^{-1}:{\mathcal{E}}_i \longrightarrow {\mathcal{E}}_0$.

Notice that one must distinguish two classes of embeddings: when the scalar
product (defined by the Killing form) in the root space $P_0$ is invariant
with respect to $\phi$ and when it is not $\phi$-invariant. The first
embedding is called "metric" , the second -- "nonmetric".

\begin{Def}
A root system $\Delta $ ''splinters'' as $(\Delta _{1},\Delta _{2})$ if
there are two embeddings $\phi _{1}:\Delta _{1}\hookrightarrow \Delta $ and $%
\phi _{2}:\Delta _{2}\hookrightarrow \Delta $ where (a) $\Delta $ is the
disjoint union of the images of $\phi _{1}$ and $\phi _{2}$ and (b) neither
the rank of $\Delta _{1}$ nor the rank of $\Delta _{2}$ exceeds the rank of $%
\Delta $.
\end{Def}

It is equivalent to say that $(\Delta_1,\Delta_2)$ is a "splint'' of $\Delta$
and we shall denote this by $\Delta \approx (\Delta_1,\Delta_2)$. Each
component $\Delta_1$ and $\Delta_2$ is a "stem'' of the splint $%
(\Delta_1,\Delta_2)$.

To study relations between injection fan technique and splint let us
consider the case when one of the stems $\Delta _{1}=\Delta _{\frak{a}}$
is a root subsystem.

Splint $\Delta \approx (\Delta _{1},\Delta _{2})$ is called ''injective'' if
$\Delta _{1}=\Delta _{\frak{a}}$, is a root subsystem
in $\Delta $ corresponding to a regular reductive subalgebra $\frak{a}%
\hookrightarrow \frak{g}$.

In case of injective splint the second stem $\Delta _{\frak{s}}:=\Delta
_{2}=\Delta \setminus \Delta _{\frak{a}}$ can be translated into a product (%
\ref{product}) and it defines an injection fan $\Gamma _{\frak{a}%
\hookrightarrow \frak{g}}$. Denote by $\Delta_{\mathfrak{s}0}$ the coimage of the second embedding $\phi:\Delta_{\mathfrak{s}0}\to \Delta_{\mathfrak{g}}$.
The following conjecture follows.

\begin{Cnj}
Each injective splint $\Delta \approx (\Delta _{\frak{a}},\Delta _{\frak{s}})
$ defines an injection fan with the carrier $\left\{ \xi \right\} _{\frak{a}%
\rightarrow \frak{g}}$ fixed by the product
\begin{equation}
\prod_{\beta \in \Delta _{\frak{s}}^{+}}\left( 1-e^{-\beta }\right)
=-\sum_{\gamma \in P}s(\gamma )e^{-\gamma }\quad   \label{splint product}
\end{equation}
\end{Cnj}

In case of injective splint we say that subalgebra $\frak{a}\hookrightarrow \frak{g}$
splinters $\Delta $ (and call $\frak{a}$ the ''splinting subalgebra'' of $%
\frak{g}$). In \cite{richter2008splints} splints are classified (see Appendix there)
and the first three types of them are injective.

\section{How stems define multiplicity functions}

\label{sec:stems and multiplicity functions}

In this Section we study properties of injective splints $\Delta
\approx (\Delta _{\frak{a}},\Delta _{\frak{s}})$. It will be
demonstrated that under certain conditions to find branching
coefficients for a splinting injection $\frak{a}\hookrightarrow
\frak{g}$ means to find weight multiplicities of an irreducible
$\frak{s}$-module $L_{\frak{s}}^{\nu }$ with fixed highest weight
$\nu $. Notice that $\frak{s}$ must not be a subalgebra of
$\frak{g}$.

Let us return to relation (\ref{singular main}) and multiply both sides by $%
R_{\frak{a}}$:
\begin{equation}
\frac{1}{\prod_{\beta \in \Delta _{\frak{s}}^{+}}(1-e^{-\beta })}\Psi _{%
\frak{g}}^{\left( \mu \right) }=\sum_{\nu \in P_{\frak{a}}^{+}}b_{\nu
}^{(\mu )}\Psi _{\frak{a}}^{\left( \nu \right) }.
\label{singular main-2}
\end{equation}
Here the first factor in the l.h.s. is the inverse of the fan $\Gamma _{%
\frak{a}\rightarrow \frak{g}}$. Consider the highest weight module $L_{\frak{%
s}}^{\nu }$. The embedding $\phi :\Delta _{\frak{s}\,0}\longrightarrow \Delta
_{\frak{g}}$ sends the singular element $\Psi _{\frak{s}}^{\left( \nu
\right) }$ into $\Psi _{\frak{g}}^{\left( \mu \right) }$. Applying the
inverse morphism $\phi ^{-1}$ to the product $\left( \prod_{\beta \in \Delta
_{\frak{s}}^{+}}(1-e^{-\beta })\right) ^{-1}\phi \left( \Psi _{\frak{s}%
}^{\left( \nu \right) }\right) $ one gets the character of the module $L_{%
\frak{s}}^{\nu }$,
\begin{equation}
\phi ^{-1}\left( \frac{1}{\prod_{\beta \in \Delta _{\frak{s}%
}^{+}}(1-e^{-\beta })}\phi \left( \Psi _{\frak{s}}^{\left( \nu \right)
}\right) \right) =\frac{1}{\prod_{\beta \in \Delta _{\frak{s}0
}^{+}}(1-e^{-\beta })}\Psi _{\frak{s}}^{\left( \nu \right) }=\mathrm{ch}%
\left( L_{\frak{s}}^{\nu }\right) .  \label{inverse for stem}
\end{equation}

Our task is to find the singular element $\Psi _{\frak{s}}^{\left(
\xi \right) }$  for the module $L_{\frak{s}}^{\xi }$ as a
component in $\Psi _{\frak{g}}^{\left( \mu \right) }$ and to prove
that $L_{\frak{s}}^{\xi }$ is uniquely defined by
$L_{\frak{g}}^{\mu }$ and that the branching coefficients $b_{\nu
}^{(\mu )}$ in the r.h.s. of (\ref {singular main-2}) coincide
with multiplicities $m_{\zeta }^{\left( \xi \right) }$ of the
corresponding weights in $\mathcal{N}_{\frak{s}}^{\xi }$ .

For a highest weight irreducible module $L_{\frak{g}}^{\mu }$ the singular
element $\Psi _{\frak{g}}^{\left( \mu \right) }$ is an element of $\mathcal{E%
}$ corresponding to the shifted Weyl-orbit of the weight $\left( \mu +\rho
\right) \in P^{+}$ with the sign function $\epsilon \left( w\right) $. It is
convenient to use also unshifted singular elements
\begin{equation}
\Phi ^{\left( \mu \right) }:=\Psi ^{\left( \mu \right) }e^{\rho }.
\label{definition Phi}
\end{equation}
In these terms the relation (\ref{singular main-2}) looks like
\begin{equation}
\frac{e^{\rho _{\frak{g}}-\rho _{\frak{a}}}}{\prod_{\beta \in \Delta _{\frak{%
s}}^{+}}(1-e^{-\beta })}\Phi _{\frak{g}}^{\left( \mu \right) }=\sum_{\nu \in
P_{\frak{a}}^{+}}b_{\nu }^{(\mu )}\Phi _{\frak{a}}^{\left( \nu \right) }.
\label{singular main-3}
\end{equation}
The orbit related to $\Phi _{\frak{g}}^{\left( \mu \right) }$ is completely
defined by the set of edges $\left\{ \lambda _{i}\right\} _{i=1,\dots ,r}$
adjusted to the end of the highest weight vector $\mu +\rho $. For $\mu
=\sum m_{i}\omega _{i}$ these edges are
\begin{equation}
\lambda _{i}=-\left( m_{i}+1\right) \alpha _{i},\quad i=1,\dots ,r.
\label{edge}
\end{equation}
Each formal exponent $e^{\mu +\rho +\lambda _{i}}$ in $\Phi _{\frak{g}%
}^{\left( \mu \right) }$ bears the sign coefficient $\epsilon =(-)$. The
defining property of $\Phi _{\frak{g}}^{\left( \mu \right) }$ is as follows.
Consider any pair of edges $\lambda _{i},\lambda _{j}$ and the corresponding
weights $\mu +\rho $, $\mu +\rho +\lambda _{i}$ and $\mu +\rho +\lambda _{j}$.
Apply the reflection $s_{\alpha _{i}}$ (or $s_{\alpha _{j}}$),
\begin{equation}
s_{\alpha _{i}}\circ \left\{
\begin{array}{l}
\left( \mu +\rho \right)  \\
\left( \mu +\rho +\lambda _{i}\right)  \\
\left( \mu +\rho +\lambda _{j}\right)
\end{array}
\right. =\left\{
\begin{array}{l}
\left( \mu +\rho +\lambda _{i}\right)  \\
\left( \mu +\rho \right)  \\
\left( \mu +\rho +\lambda _{i}-(m_{j}+1)s_{\alpha _{i}}\circ \alpha
_{j}\right)
\end{array}
\right.   \label{reflected triple}
\end{equation}

\begin{Prop}
The edge $\lambda _{i,j}$ of $\Phi _{\frak{g}}^{\left( \mu \right) }$
starting at the weight $\left( \mu +\rho +\lambda _{i}\right) $ along the
root $-s_{\alpha _{i}}\circ \alpha _{j}$ has the same length in $(s_{\alpha
_{i}}\circ \alpha _{j})$ as $\lambda _{j}$ has in $\alpha _{j}$. (The same
is true for the edge $\lambda _{j,i}$, its length in $(s_{\alpha _{j}}\circ
\alpha _{i})$ is equal to the length of $\lambda _{i}$ in $\alpha _{i}$.)
\label{diagram property}
\end{Prop}

In $\Phi _{\frak{g}}^{\left( \mu \right) }$ the elements $e^{\left( \mu
+\rho +\lambda _{i}-(m_{j}+1)s_{\alpha _{i}}\circ \alpha _{j}\right) }$ and $%
e^{\left( \mu +\rho +\lambda _{j}-(m_{i}+1)s_{\alpha _{j}}\circ \alpha
_{i}\right) }$ have the sign coefficient  $\epsilon =(+)$.

Remember that only three types of splints are
injective and thus are naturally connected with branching. Below we
reproduce the part of the splints table from \cite{richter2008splints} corresponding to
injective splints:
\begin{equation}
\label{eq:1}
\begin{array}{cc||c|c}
\hbox{type} & \hspace{0.25in}\Delta \hspace{0.25in} & \hspace{0.25in}\Delta
_{\frak{a}}\hspace{0.25in} & \hspace{0.25in}\Delta _{\frak{s}}\hspace{0.25in}
\\ \hline\hline
\hbox{(i)} & G_{2} & A_{2} & A_{2} \\
& F_{4} & D_{4} & D_{4} \\ \hline
\hbox{(ii)} & B_{r}(r\geq 2) & D_{r} & \oplus ^{r}A_{1} \\ \qquad
\qquad (*)& C_{r}(r\geq 3) & \oplus ^{r}A_{1} &  D_{r} \\ \hline
\hbox{(iii)} & A_{r}(r\geq 2) & A_{r-1}\oplus u\left( 1\right)  & \oplus
^{r}A_{1} \\
& B_{2} & A_{1}\oplus u\left( 1\right)  & A_{2}
\end{array}
\end{equation}

Each row in the table gives a splint $(\Delta _{\frak{a}},\Delta _{\frak{s}})
$ of the simple root system $\Delta $. In the first two types both $\Delta _{%
\frak{a}}$ and $\Delta _{\frak{s}}$ are embedded metrically. Stems in the
first type splints are equivalent and in the second are not. In the third
type splints only $\Delta _{\frak{a}}$ is embedded metrically. The summands $%
u\left( 1\right) $ are added to keep $r_{\frak{a}}=r$. This does not change
the principle properties of branching but makes it possible to use the
multiplicities of $\frak{s}$ -modules without further projecting their
weights. The second injective splint of type (ii) (marked by a star) does not generate a unique auxiliary $\frak{s}$ -module and in this case branching is related to splint in a more complicated form. We will not study this case here.

Splints induce a decomposition of the set $S=S_{\frak{c}}\cup S_{\frak{d}}$
with $S_{\frak{c}}=S\cap S_{\frak{a}}$ and $S_{\frak{d}}=S\cap S_{\frak{s}}$.
It is easy to check that for any injective splint the subset $S_{%
\frak{d}}$ is nonempty. It follows that in the set $\left\{ \lambda
_{i}\right\} _{i=1,\dots ,r}$ one can always find simple roots $\beta
_{k}\in \Delta _{\frak{s}}$ and that the orbit corresponding to $\Phi _{%
\frak{g}}^{\left( \mu \right) }$ contains the edges
\begin{equation}
\lambda _{k}=-\left( m_{k}+1\right) \beta _{k}  \label{beta edge}
\end{equation}
attached to the weight $\mu +\rho $. As far as $\Delta _{\frak{a}}$ is a
root system and for any pair of simple roots from $S_{\frak{c}}$ the
property \ref{diagram property} is fulfilled, the element $\Phi _{\frak{g}%
}^{\left( \mu \right) }$ being a singular element for a set of $\frak{a}$%
-modules. Consider $\beta _{l}\in \Delta _{\frak{s}}$ whose coimage in $%
\Delta _{\frak{s}0}$ is simple. In Appendix it is shown that for any such $%
\beta _{l}$ there exists a root $\alpha _{l}\in S_{\frak{c}}$ such that
$\beta _{l}=\alpha _{l}+\beta _{k}$. It is easily seen that the
corresponding edge intersects the boundary plane of the fundamental chamber $%
\bar{C_{\frak{a}}}$ orthogonal to the root $\alpha _{l}$,
\begin{equation}
s_{\alpha _{l}}\left( \mu +\rho -p\beta _{l}\right) =s_{\alpha _{l}}\left(
\mu +\rho \right) -ps_{\alpha _{l}}\beta _{l}=\mu +\rho -p\beta _{l},
\label{intersection}
\end{equation}
\begin{equation}
\mu +\rho -s_{\alpha _{l}}\left( \mu +\rho \right) =\left( m_{l}+1\right)
\alpha _{l}=\left( m_{l}+1\right) \beta _{l}-\left( m_{l}+1\right) \beta
_{k}=p\beta _{l}-ps_{\alpha _{l}}\beta _{l}.
\label{intersection-2}
\end{equation}
It follows that $p=\left( m_{l}+1\right) $ and $s_{\alpha _{l}}\beta
_{l}=\beta _{k}$. Now apply the operator $s_{\beta _{k}}$ and find that the
edge along the root $s_{\beta _{k}}\alpha _{l}$ attached at the weight $%
s_{\beta _{k}}(\mu +\rho )$ is also equal to $-ps_{\beta _{k}}\alpha _{l}$.
This means that for the triple of roots $\beta _{k},\beta _{l}$ and $%
s_{\beta _{k}}\alpha _{l}$ in $\Delta _{\frak{s}}$ the edges $\lambda
_{k}=-\left( m_{k}+1\right) \beta _{k}$, $\lambda _{l}=-\left(
m_{l}+1\right) \beta _{l}$ and $\lambda _{kl}=-\left( m_{l}+1\right)
s_{\beta _{k}}\alpha _{l}$ demonstrate the property \ref{diagram property}.
One can continue this procedure further in the 2-dimensional subspace fixed
by the roots $\beta _{k}$ and $\beta _{l}$ and find the set of formal
exponents that being supplied with the corresponding sign factors compose
the coimage of the singular element of a module for the subalgebra in $\frak{%
s}$ (this subalgebra has rank $r=2$).

The same can be proven for any positive root $\beta _{l}\in \Delta $ that is
simple in $\Delta _{\frak{s}0}$ and correspondingly for any $r=2$ subalgebra
in $\frak{s}$. The latter means that to ''find'' a singular element of $%
\frak{s}$-module in $\Phi _{\frak{g}}^{\left( \mu \right) }$ it is necessary
to incorporate in it additional formal elements $\left\{ -e^{\mu +\rho
-\left( m_{l}+1\right) \beta _{l}}|\beta _{l}\in S_{\frak{c}}\right\}.$ This
fixes the starting edges of the diagram $\phi ^{-1}\left( \Phi _{\frak{s}}^{%
\widetilde{\mu }}\right) $. As it follows from the reconstruction procedure
the highest weight $\widetilde{\mu }$ is totally defined by the weight $\mu $%
, they have the same Dynkin numbers:
\begin{equation}
\mu =\sum m_{k}\omega _{k}\qquad \Longrightarrow \quad \widetilde{\mu }=\sum
m_{k}\widetilde{\omega }_{k} . \label{new h weight}
\end{equation}

The next step is to check whether the image  $\phi \left(\Phi _{\frak{s}}^{\left(
\widetilde{\mu}\right) }\right)  $ belongs to
$\bar{C_{\frak{a}}}$ and the set $\phi
\left( \Phi _{\frak{s}}^{\widetilde{\mu }}\right) \setminus \Phi _{%
\frak{g}}^{\left( \mu \right) }|_{\bar{C_{\frak{a}}}}$ corresponds
to the weights in the boundary $\bar{C_{\frak{a}}}$ (including the
subset
$\left\{ -e^{\mu +\rho -\left( m_{l}+1\right) \beta _{l}}|\beta _{l}\in S_{%
\frak{c}}\right\} $). Provided this condition is fulfiled let us return
to relation (\ref{singular main-3}). One can add to $%
\Phi _{\frak{g}}^{\left( \mu \right) }$ pairs of formal elements
constructed above  with the opposite signs: $\epsilon \left(
w\right) |_{w\in W_{\frak{s}}}$ and $-\epsilon \left( w\right)
|_{w\in W_{\frak{s}}}$. Attribute the signs $\epsilon \left(
w\right) |_{w\in W_{\frak{s}}}$ to the elements whose weights we
shall attribute to $\bar{C_{\frak{a}}}$. The same elements with
the opposite signs are to be referred to the neighboring Weyl
chambers of $\bar{C_{\frak{a}}^{(l)}}$ (the latter are connected
with the main one via simple reflections $s_{\alpha _{l}}$ so the
signs $-\epsilon \left( w\right) |_{w\in W_{\frak{s}}}$ are
natural for them). In fact one can repeat the procedure and find
additional singular weights in any Weyl chamber
$\bar{C_{\frak{a}}^{(m)}}$ and in them additional singular weights
always have the signs opposite to that in their
nearest neighbors. Thus without changing the element $\Phi _{\frak{g}%
}^{\left( \mu \right) }$ one can present it as a sum
\begin{equation}
\Phi _{\frak{g}}^{\left( \mu \right) }=\sum_{w\in
W_{\frak{a}}}\epsilon \left( w \right) w\circ \left( e^{\rho
_{\frak{a}}}\Psi ^{\widetilde{\mu }+\rho _{\frak{s}}}\right)
\label{singular final}
\end{equation}
where the weight $\widetilde{\mu }=\sum m_{k}\omega
_{\frak{s}}^{k}$ was defined above. As far as the second condition
is fulfilled (i.e. $\phi \left(\Phi _{\frak{s}}^{\left(
\widetilde{\mu}\right) }\right) \subset \bar{C_{\frak{a}}}$) the
decomposition (\ref{singular final}) provides the
possibility to apply the factor $\left( \prod_{\beta \in \Delta _{\frak{s}%
}^{+}}(1-e^{-\beta })\right) ^{-1}$ to each summand of the singular element $%
\Phi _{\frak{g}}^{\left( \mu \right) }$ separately because the sets of
weights from different Weyl summands do not intersect. Taking into account
the isomorphism $\phi $ one can see that in the main Weyl chamber $\bar{C_{%
\frak{a}}}$ the set of weights generated by the factor $\left( \prod_{\beta
\in \Delta _{\frak{s}}^{+}}(1-e^{-\beta })\right) ^{-1}$ is isomorphic to
the weight diagram $\mathcal{N}_{\frak{s}}^{\widetilde{\mu }}$ of the $\frak{%
s}$-module $L_{\frak{s}}^{\widetilde{\mu }}$. Now one can restrict
relation (\ref{singular main-3}) to $\bar{C_{\frak{a}}}$ and obtain the main
result:
\begin{Prop}
\begin{equation}
\frac{e^{\rho _{\frak{g}}}}{\prod_{\beta \in \Delta _{\frak{s}%
}^{+}}(1-e^{-\beta })}\left( \Psi ^{\widetilde{\mu }+\rho _{\frak{s}%
}}\right) =\sum_{\widetilde{\nu }\in \mathcal{N}_{\frak{s}}^{\widetilde{\mu }%
}}M_{\left( \frak{s}\right) \widetilde{\nu }}^{\widetilde{\mu }}e^{\left(
\mu -\phi \left( \widetilde{\mu }-\widetilde{\nu }\right) \right)
}=\sum_{\nu \in P_{\frak{a}}^{++}}b_{\nu }^{(\mu )}e^{\nu }.
\label{singular main-4}
\end{equation}
Any weight with nonzero multiplicity in the r. h. s. is equal to one of the
highest weights in the decomposition. The multiplicity $M_{\left( \frak{s}%
\right) \widetilde{\nu }}^{\widetilde{\mu }}$ of the weight  $\widetilde{\nu
}\in \mathcal{N}_{\frak{s}}^{\widetilde{\mu }}$ defines the branching
coefficient $b_{\nu }^{(\mu )}$ for the highest weight $\nu =\left( \mu
-\phi \left( \widetilde{\mu }-\widetilde{\nu }\right) \right) $:
\[
b_{\left( \mu -\phi \left( \widetilde{\mu }-\widetilde{\nu }\right) \right)
}^{(\mu )}=M_{\left( \frak{s}\right) \widetilde{\nu }}^{\widetilde{\mu }}.
\]
\end{Prop}

\section{Examples}
\label{sec:examples}
\begin{example}
  Consider the Lie algebra $A_{2} =\bf{sl}(3)$ and branching of its irreducible module $L^{[3,2]}_{A_{2}}$ with respect to the reductive subalgebra $A_{1}\oplus u(1)$. The root system  $\Delta _{\frak{a}}= \Delta_{A_{1}\oplus u(1)}$ contains the simple root of $\alpha_1=e_1-e_2$ of $A_{2}$. The singular element $\Psi^{[3,2]}_{\frak{a}}$ is decomposed into a sum of splint images of singular elements of $A_{1}\oplus A_{1}$-modules. Branching coefficients $b_{\nu }^{[ 3,2 ] }$ coincide with weight multiplicities of $L^{[3,2]}_{A_{1}\oplus A_{1}}$-module (see Fig. \ref{fig:a2_splint}).
  \begin{figure}[h!bt]
  \noindent\centering{
   \includegraphics[width=75mm]{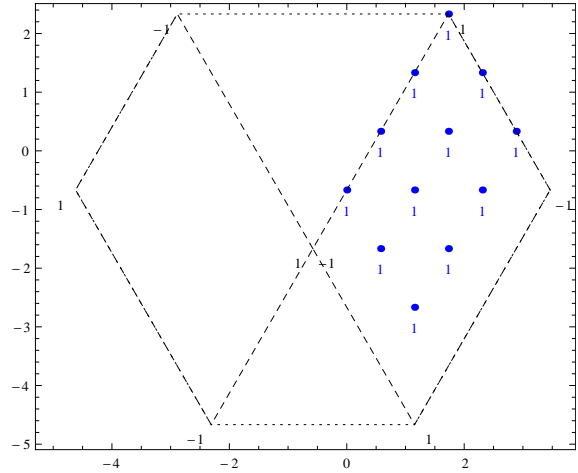}
  }
  \caption{Weyl group orbit (dotted) producing singular element of $L^{[3,2]}_{A_{2}}$ and its decomposition into the sum of splint images of singular elements of modules  $L^{[3,2]}_{A_{1}\oplus A_{1}}$ (dashed). Weight multiplicities of $L^{[3,2]}_{A_{1}\oplus A_{1}}$-module coincide with branching coefficients for the reduction $L^{[3,2]}_{A_{2}\downarrow A_{1}\oplus u(1)}$.}
\label{fig:a2_splint}
\end{figure}
\end{example}

\begin{example}
  For the Lie algebra $B_{2}= \bf{so}(5)$ branching of its irreducible module $L^{[3,2]}$ into  modules of a reductive subalgebra $A_{1}\oplus u(1)$ with the root system spanned by the first simple root $\alpha_1=e_1-e_2$ of $B_{2}$. Singular element of $\Psi^{[3,2]}_{B_{2}}$ is decomposed into the sum of splint images of singular elements of $A_{2}$-modules and branching coefficients coincide with weight multiplicities of $A_{2}$-module (see Fig. \ref{fig:b2_splint}).

  \begin{figure}[h!bt]

   \includegraphics[width=65mm]{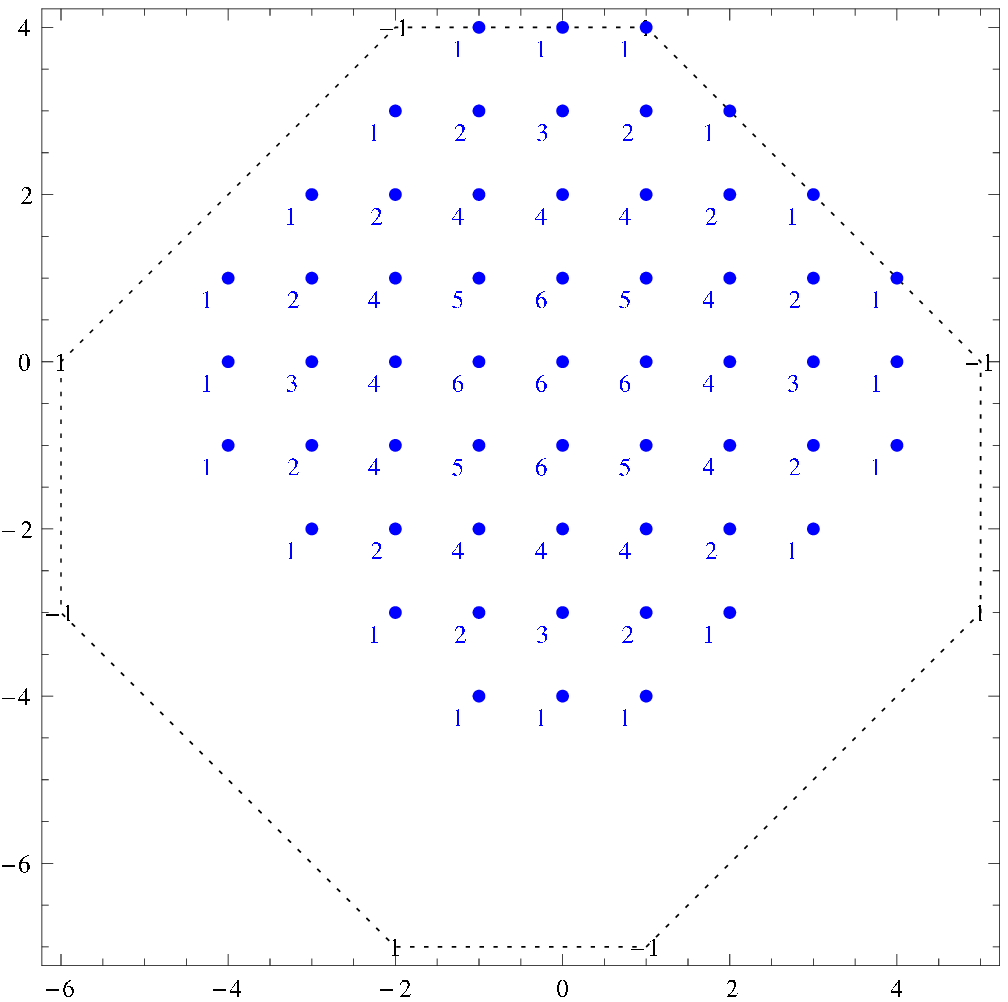}
   \includegraphics[width=65mm]{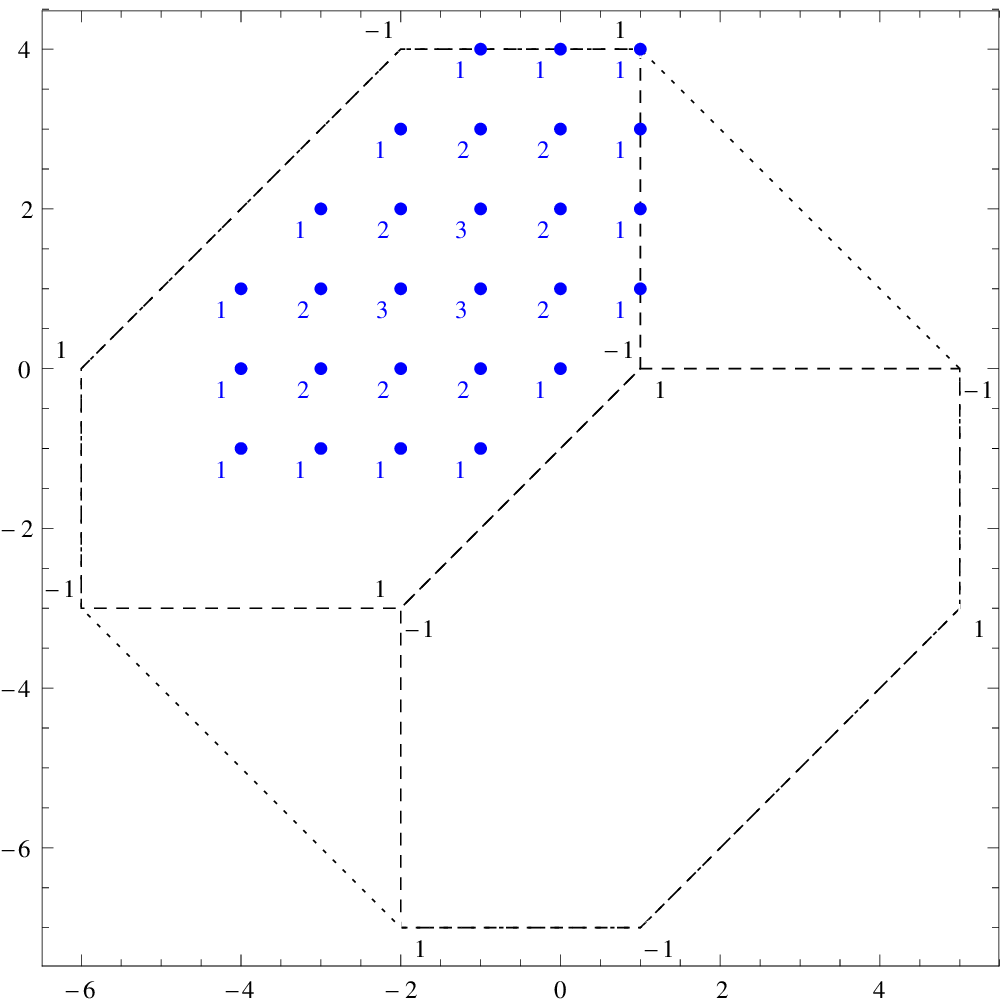}
  \caption{Weights of the $B_{2}$-module $L^{[3,2]}$ are indicated by dots in the left picture (their multiplicities are also indicated). Contour of the singular element is shown by dotted line. The right picture presents the decomposition of  $\Psi_{B_{2}}(L^{[3,2]}_{B_{2}})$-singular element into the sum of splint images of singular elements $\Psi_{A_{2}}(L^{[3,2]})$ (dashed). Weight multiplicities of $L^{[3,2]}_{A_{2}}$-module coincide with branching coefficients for the reduction $L^{[3,2]}_{B_{2}\downarrow A_{1}\oplus u(1)}$.}

 \label{fig:b2_splint}
\end{figure}
\end{example}

\vspace{10mm}
\begin{example}
   Lie algebra $G_{2}$ has a regular subalgebra $A_{2}$ with root system $\Delta_{\mathfrak{a}}=\Delta_{A_{2}}$ containing the $G_{2}$ long roots. Consider branching of an irreducible module $L_{G_{2}}^{(3,2)}$ into the $A_{2}$-modules. Singular element $\Psi_{G_{2}}(L^{[3,2]})$ is decomposed into the sum of splint images of singular elements $\Psi_{A_{2}}(L^{[3,2]})$ and the corresponding branching coefficients coincide with weight multiplicities of $L^{[3,2]}_{A_{2}}$-module (see Fig. \ref{fig:g2_splint}).

  \begin{figure}[h!bt]
  \noindent\centering{
   \includegraphics[width=120mm]{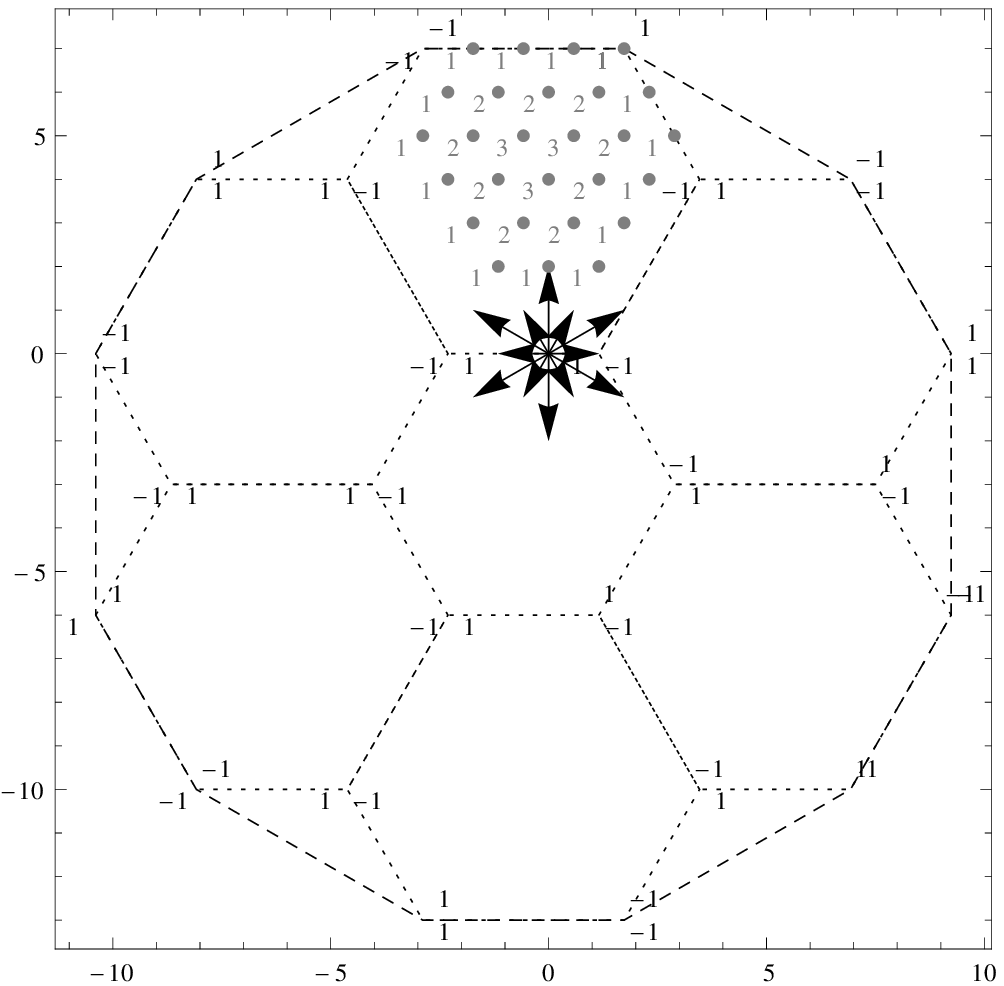}
  }

  \caption{Weyl group orbit (dotted) for the singular element $\Psi_{G_{2}}(L^{[3,2]})$ and its decomposition into the sum of splint images of singular elements of $A_{2}$-modules (dashed). Weight multiplicities of $L^{[3,2]}_{A_{2}}$-module coincide with branching coefficients for the reduction $L^{[3,2]}_{G_{2}\downarrow A_{2}}$.}

 \label{fig:g2_splint}
\end{figure}

\end{example}

\section{Conclusions}

\label{sec:conclusions}It is explicitly demonstrated that splint
presents a very effective tool to find branching coefficients. In
particular the injective splints that have the property $\phi
\left( \Phi _{\frak{s}}^{\left( 0 \right)}\right) \subset
\bar{C_{\frak{a}}^{(0)}}$ provide the possibility to reduce
branching rules calculations for the highest weight modules to a
determination of weight multiplicities for a module with the same
Dynkin labels referred to the Lie algebra $\mathfrak{s}$. This
algebra $\mathfrak{s}$ must not be a subalgebra in the initial
$\mathfrak{g}$, it has the same rank $r_{\mathfrak{s}}=r$, but
obviously is less "complicated" than $\mathfrak{g}$ -- only a
subset of the initial root system is involved in the second stem
$\Delta_{\mathfrak{s}}$.

It is significant that for the injections $D_{r}\hookrightarrow
B_{r}$ and $A_{r-1}\oplus u\left( 1\right) \hookrightarrow A_{r}$
the splint technique shows transparently Gelfand-Tzeytlin rules
for branching:  the reduction is multiplicity free (all nonzero
branching coefficients are equal to 1). Here it is an immediate
consequence of the structure of the second stem being a direct sum
of $A_{1}$ algebras and the fact that the corresponding module
$L_{\frak{s}}^{\mu }$ is irreducible.

\section*{Acknowledgements}
\label{sec:acknowledgements} The authors are grateful to Prof.
David Richter for his important notes. The work was supported in
part by the RFFI grant N 09-01-00504. A.A.N. thanks the Chebyshev
Laboratory (Department of Mathematics and Mechanics,
Saint-Petersburg State University) for support under the grant
11.G34.31.0026 of the Government of Russian Federation.

\bibliography{bibliography}{}

\providecommand{\href}[2]{#2}\begingroup\raggedright\begin{thebibliography}{1}

\bibitem{richter2008splints}
D.~Richter, ``Splints of classical root systems,'' {\em Arxiv preprint
  arXiv:0807.0640} (2008)  , \href{http://arxiv.org/abs/0807.0640}{{\tt
  0807.0640}}.

\bibitem{lyakhovsky1996rra}
V.~Lyakhovsky, S.~Melnikov, {\em et al.}, ``{Recursion relations and branching
  rules for simple Lie algebras},'' {\em Journal of Physics A-Mathematical and
  General} {\bf 29} (1996) no.~5, 1075--1088,
  \href{http://arxiv.org/abs/q-alg/9505006}{{\tt q-alg/9505006}}.

\bibitem{2010arXiv1007.0318L}
V.~{Lyakhovsky} and A.~{Nazarov}, ``{Recursive algorithm and branching for
  nonmaximal embeddings},''
  \href{http://dx.doi.org/10.1088/1751-8113/44/7/075205}{{\em Journal of
  Physics A: Mathematical and Theoretical} {\bf 44} (2011) no.~7, 075205},
  \href{http://arxiv.org/abs/1007.0318}{{\tt arXiv:1007.0318 [math.RT]}}.

\bibitem{ilyin812pbc}
M.~Ilyin, P.~Kulish, and V.~Lyakhovsky, ``{On a property of branching
  coefficients for affine Lie algebras},'' {\em Algebra i Analiz} {\bf 21}
  (2009)  2, \href{http://arxiv.org/abs/0812.2124}{{\tt arXiv:0812.2124
  [math.RT]}}.

\bibitem{humphreys1997introduction}
J.~Humphreys, {\em {Introduction to Lie algebras and representation theory}}.
\newblock Springer, 1997.

\end{thebibliography}\endgroup
\bibliographystyle{utphys}


\section*{Appendix}


Let us demonstrate that for injective splints of classical Lie
algebras the following property is valid:
\begin{Prop}
Let $\Delta \approx (\Delta _{\frak{a}},\Delta _{\frak{s}})$ be an
injective
splint with the decomposition of simple roots $S=S_{\frak{c}}\cup S_{\frak{d}%
}$ with $S_{\frak{c}}=S\cap S_{\frak{a}}$ and $S_{\frak{d}}=S\cap S_{\frak{s}%
}$.
For any simple root $\beta \in S_{\frak{s}}$ there exists the
pair of
roots ( $\alpha $ ,$\beta ^{\prime }$) with $\alpha \in $ $S_{\frak{c}%
},\beta ^{\prime }\in S_{\frak{s}}$ such that $\alpha =\beta
-\beta ^{\prime }$.
\end{Prop}

\begin{itemize}
\item
Type 1. $\Delta _{G_{2}}\approx (\Delta _{A_{2}},\Delta
_{A_{2}}).$

 \begin{figure}[h!bt]
  \noindent\centering{
   \includegraphics[width=100mm]{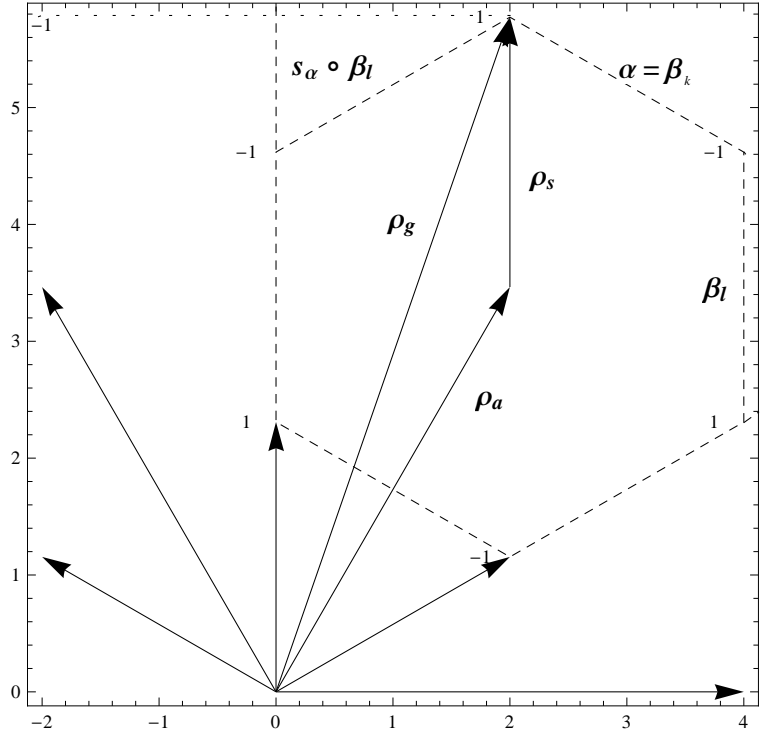}
  }
  \caption{Positive roots of $G_{2}$ and formation of singular element $\Phi^{(0)}_{\mathfrak{s}}$ in the main Weyl chamber of $\mathfrak{a}=A_{2}$.}
\end{figure}
Here both stems are metric and the corresponding root systems are
equivalent. In Figure 4 a part of the singular element $\Phi
_{G_{2}}^{\left( 0\right) }$ is presented. The boundaries of $\bar{C_{\frak{a%
}}}$ are the dashed lines starting at the center of the singular
element. It
contains the edge $\lambda _{2}=-\alpha _{2}=-\beta _{2}$ and the roots $%
-\beta _{1}=-s_{\alpha _{2}}\circ \beta _{3}$ and $ -\beta _{3}$
($\beta _{3}$ is indicated as $\beta _{l}$). For the root $\beta
_{1}$ the necessary pair
is  $(\alpha _{1}, \beta _{2})$: $\alpha _{1}=\beta _{1}-\beta _{2}$. The $%
\lambda _{2,3}^{\frak{s}}=\beta _{3}$ edge is equal to $\lambda _{1}^{\frak{s%
}}=\beta _{1}=s_{\alpha _{2}}\circ \beta _{3}$ and $m_{1}$ index
is aquired by the $\frak{s}$-module that also inherit the second
index $m_{2}$. In this particular
case they are $m_{1}=m_{2}=0$. The general case with the initial module $%
L^{\mu }$ and $\mu =m_{1}\omega _{1}+m_{2}\omega _{2}$ can be
treated in the same way: one finds an edge $\lambda _{2}=-\left(
m_{2}+1\right) \beta _{2}$ and put $\lambda
_{1}^{\frak{s}}=-\left( m_{1}+1\right) \beta _{1}$, its end
belongs to the boundary $\bar{C_{\frak{a}}}$ . The reflection
$s_{\beta
_{2}} $ sends $\beta _{1}$ to $\beta _{3}$ and the corresponding edge
$\lambda _{2,3}^{\frak{s}}=-\left( m_{1}+1\right) \beta _{3}$ has
the length $\left( m_{1}+1\right) $. Now consider $\lambda _{1}^{\frak{s}}$ (or
$\lambda _{2,3}^{\frak{s}} $) and $\lambda _{1,3}^{\frak{s}}$ (or
$\lambda _{2,3,1}^{\frak{s}}$) edges to find that they belong to
the boundary $\bar{C_{\frak{a}}}$ and the Weyl symmetry predicts
that $\lambda _{1,3}^{\frak{s}}=-\left( m_{2}+1\right)
\beta _{3}$ ($\lambda _{2,3,1}^{\frak{s}}=-\left( m_{2}+1\right) \beta _{1}$%
) . Finally the edge $\lambda _{1,3,2}^{\frak{s}}=-\left(
m_{1}+1\right) \beta _{2}$ closes the polytope. Its vertices
correspond to weights of the singular
element $\Phi _{\frak{s}}^{\left( \widetilde{\mu }\right) }=\sum_{w\in W_{%
\frak{s}}}\varepsilon \left( w\right) e^{w\circ \left(
\widetilde{\mu }+\rho
_{\frak{s}}\right) }$ of the module $L_{\frak{s}}^{\left( \widetilde{\mu }%
\right) }$ with $\widetilde{\mu }=m_{1}\widetilde{\omega }_{1}+m_{2}%
\widetilde{\omega }_{2}$ . Notice that in this case the sign
factors can be obtained directly in the initial weight system as
far as the stem is metric.
\item
Type 1. $\Delta _{F_{4}}\approx (\Delta _{D_{4}},\Delta
_{D_{4}}).$

Both stems are metric here and the corresponding root systems are
equivalent. The system $\Delta _{D_{4}}$ of the subalgebra
$\frak{a=}D_{4}$ is formed by the set $\left\{ \pm e_{i}\pm
e_{j}\right\} _{|i,j=1,\ldots 4,\; i\neq j}.$ The simple roots
$S_{\frak{c}}$ are $\left\{ e_{2}-e_{3},e_{3}-e_{4}\right\} $ and
$S_{\frak{d}}=\left\{ e_{4},\frac{1}{2}\left(
e_{1}-e_{2}-e_{3}-e_{4}\right) \right\} $.
For a module $L^{\mu }$ with $\mu =\sum m_{k}\omega _{k}$ consider the edge $%
\lambda _{3}=-\left( m_{3}+1\right) e_{4}=-\left( m_{3}+1\right)
\beta _{3}$.
 Compose an edge $\lambda _{2}^{\frak{s}}=-\left( \widetilde{m}%
_{2}+1\right) \beta _{2}$. The necessary pair of roots is $\left(
\alpha
_{2}=e_{3}-e_{4},\beta _{3}\right) $. The intersection of $\lambda _{2}^{%
\frak{s}}$ with the $\alpha _{2}$ -boundary of
$\bar{C_{\frak{a}}}$ fixes
its length to be $\lambda _{2}^{\frak{s}}=-\left( m_{2}+1\right) \beta _{2}$%
\ and the length of the edge $\lambda _{3,2}^{\frak{s}}$ is equal
to that of
$\lambda _{2}^{\frak{s}}$. Next consider the edge $\lambda _{2}^{\frak{s}%
}=-\left( m_{2}+1\right) \beta _{2}$ and the pair $\left( \alpha
_{1}=e_{2}-e_{3},\beta _{1}=e_{2}\right) $. The length of
$\lambda _{1}^{\frak{s}}$ becomes equal to $\lambda
_{1}^{\frak{s}}=-\left( m_{1}+1\right) \beta _{1}$. Proceed further till the closure of the polytope.
The edges looking along the roots of the $\alpha_4$-type, $\alpha _{4}=\beta
_{4}=$ $\frac{1}{2}\left( e_{1}-e_{2}-e_{3}-e_{4}\right) $, are
treated similarly and finally the singular element $\Phi
_{\frak{s}}^{\left( \widetilde{\mu }\right) }=\sum_{w\in
W_{\frak{s}}}\varepsilon \left( w\right) e^{w\circ \left(
\widetilde{\mu }+\rho _{\frak{s}}\right) }$ for the module
$L_{\frak{s}}^{\left( \widetilde{\mu }\right) }$ with
$\widetilde{\mu }=\sum m_{k}\widetilde{\omega }_{k}$ is formed in
$\bar{C_{\frak{a}}}$.
\item
Type 2. $\Delta _{B_{r}}\approx (\Delta _{D_{r}},\Delta _{\oplus
^{r}A_{1}}). $

Both stems are metric. An injection is fixed by the stem $\Delta
_{D_{r}}$ simple roots $S_{\frak{a}}=\left\{
e_{1}-e_{2},e_{2}-e_{3},\ldots
,e_{r-1}-e_{r},e_{r-1}+e_{r}\right\} $ . The second stem
corresponds to a
direct sum of algebras $A_{1}$ with the simple roots $S_{\frak{s}%
}=\left\{ e_{1},e_{2},\ldots ,e_{r-1},e_{r}\right\} $. Consider the edge $%
\lambda _{r}=-\left( m_{r}+1\right) \beta _{r}$ (here $\beta
_{r}=e_{r}$) and $\lambda _{r-1}=-\left(
\widetilde{m}_{r-1}+1\right) \beta _{r-1}$
attached to it (here $\beta _{r-1}=e_{r-1}$). The corresponding pair is $%
\left( \alpha _{r-1}=e_{r-1}-e_{r},\beta _{r-1}=e_{r-1}\right) $.
 The intersection condition fixes the second edge to be $\lambda
_{r-1}=-\left( m_{r-1}+1\right) \beta _{r-1}$ , it is orthogonal
to $\beta _{r}$ so the opposite edge has the same length. The
Dynkin index $m_{r-1}$ now refers
also to the\ simple root $\beta _{r-1}$ . Next consider the obtained edge $%
\lambda _{r-1}=-\left( m_{r-1}+1\right) \beta _{r-1}$ and $\lambda
_{r-2}=-\left( \widetilde{m}_{r-2}+1\right) \beta _{r-2}$ to
fix the index $\widetilde{m}_{r-2}=m_{r-2}$ and the edge
$\lambda _{r-2}=-\left( m_{r-2}+1\right) \beta _{r-2}$ and so on
till all the pairs of edges are properly fixed. Finally in
$\bar{C}_{D_{r}}$ the element $\Phi _{\oplus ^{r}A_{1}}^{\left(
\widetilde{\mu }\right) }=\sum_{w\in W_{\oplus
^{r}A_{1}}}\varepsilon \left( w\right) e^{w\circ \left( \widetilde{\mu }+%
\frac{1}{2}\sum e_{k}\right) }$ can be constructed for the module
$L_{\oplus ^{r}A_{1}}^{\left( \widetilde{\mu }\right) }$ with
$\widetilde{\mu }=\sum m_{k}\frac{1}{2}e_{k}$ .

\item
Type 2. $\Delta _{C_{r}}\approx (\Delta _{\oplus
  ^{r}A_{1}},\Delta _{D_{r}}). $

The situation in this case is analogous to the previous one and
the additional edges are constructed similarly. However in this
case the property $\phi \left( \Phi _{\frak{s}}^{\left( 0
\right)}\right) \subset \bar{C_{\frak{a}}^{(0)}}$ is violated. The
set $\phi \left( \Phi _{\frak{s}}^{\left( 0 \right)}\right)$
contains weights in several bordering  Weyl chambers
$C_{\frak{a}}$. The decomposition (\ref{singular final}) cannot be
performed. The injective splint $\Delta _{C_{r}}\approx (\Delta
_{\oplus^{r}A_{1}},\Delta _{D_{r}}) $ does not induce the property
(\ref{singular main-4}).

\item
Type 3 $\Delta _{A_{r}}\approx (\Delta _{A_{r-1}\oplus
u_{1}},\Delta _{\oplus ^{r}A_{1}}).$

Here only the first stem is metric and it fixes the injection with
simple roots $S_{\frak{a}}=\left\{ e_{1}-e_{2},e_{2}-e_{3},\ldots
,e_{r-1}-e_{r}\right\} $. The second stem corresponding to a direct sum of $%
r$ copies of $A_{1}$ has the simple roots $S_{\frak{s}}=\left\{
e_{1}-e_{r+1},e_{2}-e_{r+1},\ldots ,e_{r}-e_{r+1}\right\} $.
Consider the edge $\lambda _{r}=-\left( m_{r}+1\right) \beta _{r}$
with $\beta
_{r}=e_{r}-e_{r+1}$ and $\lambda _{r-1}=-\left( \widetilde{m}%
_{r-1}+1\right) \beta _{r-1}$ \ with $\beta
_{r-1}=e_{r-1}-e_{r+1}$ attached to it. Then the corresponding
pair is $\left( \alpha _{r-1}=e_{r-1}-e_{r},\beta
_{r-1}=e_{r-1}-e_{r+1}\right)$. The intersection with the
boundary of $\bar{C}_{A_{r-1}}$ orthogonal to $\alpha _{r-1}$
fixes the second edge to be $\lambda _{r-1}=-\left(
m_{r-1}+1\right) \beta _{r-1}$. The Dynkin index $m_{r-1}$ is to
be used for the fundamental weight  $\omega _{r-1}.$ The\
reflection $s_{\beta _{r}}$ sends $\lambda _{r-1}=-\left(
m_{r-1}+1\right) \beta _{r-1}$ to $\lambda _{r,r-1}=-\left(
m_{r-1}+1\right) \beta _{r-1}.$ Next consider the obtained edge
$\lambda _{r-1}=-\left( m_{r-1}+1\right) \beta _{r-1}$ and
$\lambda _{r-2}=-\left( \widetilde{m}_{r-2}+1\right) \beta _{r-2}$
with $\beta _{r-2}=e_{r-2}-e_{r+1} $ to obtain the index
$\widetilde{m}_{r-2}=m_{r-2}$ and the edge $\lambda _{r-2}=-\left(
m_{r-2}+1\right) \beta _{r-2}$ and so on till all the pairs of
edges are properly fixed. Finally in $\bar{C}_{D_{r}}$ the element
$\Phi _{\oplus ^{r}A_{1}}^{\left( \widetilde{\mu }\right)
}=\sum_{w\in W_{\oplus
^{r}A_{1}}}\varepsilon \left( w\right) e^{w\circ \left( \widetilde{\mu }+%
\widetilde{\rho }\right) }$ can be constructed for the module
$L_{\oplus ^{r}A_{1}}^{\left( \widetilde{\mu }\right) }$ with
$\widetilde{\mu }=\sum m_{k}\beta _{k}$. The simplest case $\Delta
_{A_{2}}\approx (\Delta _{A_{1}\oplus u_{1}},\Delta _{A_{1}\oplus
A_{1}})$ is presented in Example 4.1 and Figure 1.
\item
Type 3 $\Delta _{B_{2}}\approx (\Delta _{A_{1}},\Delta _{A_{2}}).$

This splint is illustrated in Example 4.1 and Figure 1,
$S_{A\_1}=\left\{ e_{1}-e_{2}\right\} $, $S_{A\_2}=\left\{
e_{1},e_{2}\right\} $. The edge $\lambda _{\alpha _{2}}=\lambda
_{\beta _{2}}=-\left( m_{2}+1\right) \beta _{2}$ is followed by
$\lambda _{\beta _{1}}=-\left( \widetilde{m}_{1}+1\right) \beta
_{1}$. Consider the pair $\left( \alpha _{1}=e_{1}-e_{2},\beta
_{1}=e_{1}\right)$. The end of the edge $\lambda _{\beta _{1}}$
must indicate a weight invariant under the reflection $s_{\alpha
_{1}}$.  Its length is thus fixed: $\lambda _{\beta _{1}}=-\left(
m_{1}+1\right) \beta _{1}$. In the coimage of the second stem,
that is in the root system $\Delta_{A_{2}}$, the reflection
$s_{\beta _{2}}$ sends $\lambda _{\beta _{1}}=-\left(
m_{1}+1\right) \beta
_{1}$ to $\lambda _{2,3}$, thus the latter edge has the same length in $%
\beta _{3}=e_{1}+e_{3}$, we have $\lambda _{2,3}=-\left(
m_{1}+1\right)
\beta _{3}$ with $\beta _{3}=e_{1}+e_{3}$. The irreducible $\frak{s}$%
-module has the highest weight $\widetilde{\mu }=m_{1}\widetilde{\omega }%
_{1}+m_{2}\widetilde{\omega }_{2}$. In Figure 1 we see the details
of these relations in a particular case where $L_{B_{2}}^{\left[
3,2\right] }$ is reduced to a subalgebra $A_{1}\oplus u\left(
1\right) $ and the
corresponding highest weights (with their multiplicities) form the diagram $%
\mathcal{N}_{A_2}^{\left[ 3,2\right] }$.
\end{itemize}


\end{document}